\documentclass[a4paper,10pt,reqno]{amsproc}

\usepackage[textwidth=11.4cm,textheight=18.1cm]{geometry}
\usepackage{amsfonts,amsmath,amssymb,amsthm}
\usepackage{graphicx}
\usepackage{fixmath}
\usepackage[hidelinks]{hyperref}
\usepackage[english]{babel}

\usepackage{setspace}
\usepackage{subfig}
\usepackage{url}
\usepackage{booktabs}
\usepackage{ifthen}
\usepackage{wrapfig}

\usepackage{enumerate}

\usepackage{tikz}
\usetikzlibrary{calc}
\usetikzlibrary{decorations.pathreplacing}

\usepackage[T2A,T1]{fontenc}

\newtheorem{thm}{Theorem}[section]

\newtheorem{defn}[thm]{Definition}

\newtheorem{q}{Question}

\numberwithin{equation}{section}

\begin{document}
 
\thispagestyle{plain}
 
\title{A counterexample to a conjecture of Ghosh}

\author{Hung Hua}
\address{Hung Hua   {\texttt hhua@student.clayton.edu}}

\author{Elliot Krop}
\address{Elliot Krop   {\texttt elliotkrop@clayton.edu}}

\author{Christopher Raridan}
\address{Christopher Raridan   {\texttt christopherraridan@clayton.edu}}

\date{\today}

\begin {abstract}
We answer two questions of Shamik Ghosh in the negative. We show that there exists a lobster tree of diameter less than $6$ which accepts no $\alpha$-labeling with two central vertices labeled by the critical number and the maximum vertex label. We also show a simple example of a tree of diameter $4$, with an even degree central vertex which does not accept a maximum label in any graceful labeling.
\\[\baselineskip] 2010 Mathematics Subject
      Classification: 05C78
\\[\baselineskip]
\end {abstract}

\maketitle

 \section{Introduction}
 For basic graph theoretic notation and definitions we refer to West~\cite{West}. The well-known graceful tree conjecture states that for any tree $T$ on $n$ vertices, there exists an injective vertex-labeling from $\{0,\dots, n-1\}$, so that the set of  edge-weights, defined as absolute difference on the labels on incident vertices for each edge, is $\{1,\dots, n-1\}$.

We call a labeling $f$ \emph{bipartite}, if  there exists and integer $c$ such that for any edge $uv$, either $f(u)\leq c < f(v)$ or $f(v)\leq c < f(u)$. A bipartite labeling that is graceful is called an $\alpha$\emph{-labeling}. We call the number $c$, the \emph{critical number}.

For any tree $T$, let $P$ be a longest path in $T$ and call $T$ a $k$-distant tree if all of its vertices are a distance at most $k$ from $P$. At this time, the conjecture is unknown for any trees with tree distance at least $2$. For $2$-distant trees, or \emph{lobsters}, the problem is known as Bermond's conjecture\cite{Bermond}.
 
Recently, in the pursuit of obtaining graceful lobsters from smaller graceful lobster, Ghosh \cite{Ghosh} asked two questions, which if answered in the affirmative, could have led to the solution of Bermond's conjecture. We address both of these questions:
\begin{q}\cite{Ghosh}
Does there exist an $\alpha$-labeling of a lobster of diameter at most $5$ such that the central vertices are labeled by the critical number and the maximum label?
\end{q}

\begin{q}\cite{Ghosh}
For any tree $T$ of diameter $4$, with central vertex $v$ of even degree, does there exist a graceful labeling of $T$ such that the label on $v$ is the maximum label?
\end{q}

Note: The trees considered in Question $1$ were those with two central vertices. In the argument below we answer this as well as another related question.

\begin{defn}
A vertex $v$ in a tree $T$ is an \emph{almost central vertex}, if $v$ is adjacent to a central vertex and lies on a longest path of $T$.
\end{defn}

\begin{q}
Does there exist an $\alpha$-labeling of a lobster of diameter at most $5$ such that the central vertex and an almost central vertex is labeled by the critical number and the maximum label, respectively?
\end{q}

\section{An example}
Let $T$ be the following simple $1$-distant tree.
\begin{figure}[ht]
\begin{center}
\begin{tikzpicture}[scale=1.5]
\tikzstyle{vert}=[circle,fill=black,inner sep=3pt]
\tikzstyle{overt}=[circle,fill=black!30, inner sep=3pt]

  \node[vert, label=left:\tiny{}] (u-1) at (1,1) {};
  \node[overt, label=above:{$v_1$}] (u-2) at (2,1) {};
  \node[vert, label=above:{$v$}] (u-3) at (3,1) {};
  \node[overt, label=above:{$v_2$}] (u-4) at (4,1) {};
  \node[vert, label=right:\tiny{}] (u-5) at (5,1) {};

  \node[vert, label=left:\tiny{}] (v-1) at (1,1.5) {};

   \draw[color=black] 
   (u-1)--(u-2)-- (u-3)--(u-4) --(u-5)
   (v-1)--(u-2)
;

\end{tikzpicture}
\caption{$T$}
\label{T}
\end{center}
\end{figure}

Van Bussel \cite{VB} showed that $T$ is not $0$-centered, that is, does not accept a graceful labeling with central vertex $v$ labeled $0$. By considering the complementary labelings (vertices relabeled by $n-1$ minus their old labels), we conclude that the central vertex cannot be labeled by the maximum label, $n-1$. This observation gives a negative answer to Question $2$. 

We consider Question $3$ and again refer to $T$. If this question is to be answered in the affirmative, $v$ must be labeled by the critical number $c$. Since we are looking for an $\alpha$-labeling, and if $v$ is in the lesser labeled partite set, we must have $c=3$. There are only a few cases left to consider; either $v_1$ or $v_2$ is labeled by the maximum label $5$. Since the vertex with the maximum label must be adjacent to the vertex labeled $0$, it is easy to verify that there is no graceful labeling of $T$ whether we label $v_1$ or $v_2$ by $5$.

To answer Question $1$, we need another simple example.

\vspace{1 in}

Let $S$ be the following simple $1$-distant tree.
\begin{figure}[ht]
\begin{center}
\begin{tikzpicture}[scale=1.5]
\tikzstyle{vert}=[circle,fill=black,inner sep=3pt]
\tikzstyle{overt}=[circle,fill=black!30, inner sep=3pt]

  \node[vert, label=left:{$u_2$}] (u-1) at (1,1) {};
  \node[overt, label=above:{$v_1$}] (u-2) at (2,1) {};
  \node[vert, label=above:{$v$}] (u-3) at (3,1) {};
  \node[overt, label=above:{$v_2$}] (u-4) at (4,1) {};
  \node[vert, label=above:{$u_3$}] (u-5) at (5,1) {};
  \node[overt, label=above:{$v_3$}] (u-6) at (6,1) {};    
    
  \node[vert, label=left:{$u_1$}] (v-1) at (1,1.5) {};

   \draw[color=black] 
   (u-1)--(u-2)-- (u-3)--(u-4) --(u-5)--(u-6)
   (v-1)--(u-2)
;

\end{tikzpicture}
\caption{$S$}
\label{S}
\end{center}
\end{figure}

We notice that for an $\alpha$-labeling of $S$, the bipartitions are $\{u_1,u_2,v,u_3\}$ and $\{v_1,v_2,v_3\}$.

 Assume first that $v$ is labeled by the maximum label, $6$. In this case the critical number is $2$, so to satisfy the conditions of Question $1$, label $v_2$ by $2$. Since the vertex with the maximum label must be adjacent to the vertex labeled $0$ in any graceful labeling, $v_1$ must be labeled $0$. Since $v_3$ is in the bipartition with the lower labels, the label $1$ is forced on $v_3$. Notice that in this case, $u_3$ cannot be labeled $5$, otherwise the weight $4$ would be on two edges. Furthermore, it is easy to check that either a label of $3$ or $5$ on $u_3$ would produce two edges of the same weight in $S$.

Next, consider the case when $v_2$ is labeled by $6$. In this case, the critical number is $3$, so we label $v$ by $3$. Arguing as above, $u_3$ must be labeled $0$, which leads to labeling $u_1$ and $u_2$ by $1$ and $2$. This leaves only two possibilities. Either $v_1$ is labeled $4$ or $5$, and in either case the labeling is not graceful.

We note, as suggested by the anonymous referee, that this example is vertex minimum in the sense of Question $1$. That is, it is impossible to delete $u_1$ (for example) since this would produce $P_6$, which does accept an $\alpha$-labeling such that the central vertices are labeled by the critical number and the maximum label. Furthermore, there are $\alpha$-labelings of $S$ in which $v_2$ receives either the critical number or the maximum label.

\section{Acknowledgements}
We would like to thank Shamik Ghosh for his valuable comments. We would also like to thank the anonymous referee for the explanation of Question $1$, an argument for why the answer is no, and for providing such a helpful review.

 \bibliographystyle{plain}

\begin{thebibliography}{MMM}

\bibitem{Bermond} J.-C.~Bermond, Graceful graphs, radio antennae, and French windmills, In R.J. Wilson, editor, {\it Graph Theory and Combinatorics}, p. 18-37. Pitman Publishing Ltd., (1979).

\bibitem{Ghosh} S.~Ghosh, On certain classes of graceful lobsters, submitted, {\it arXiv:1306.2932 [math.CO]} (2013)

\bibitem{VB} F.~Van Bussel, 0-Centred and 0-ubiquitously graceful trees, {\it Discrete Mathematics} 277(1-3): 193-218 (2004)

\bibitem{West}D.B.~West, {\it Introduction to Graph Theory}, second edition, Prentice-Hall (2001).
 
 \end{thebibliography}
 
 \end{document}